\input amstex
\documentstyle{amsppt}
\refstyle{B}


\topmatter
\title
{Eigenbundles,  Quaternions, and Berry's Phase}
\endtitle

\author
Daniel Henry Gottlieb
\endauthor


\address
Math. Dept., Purdue University, West Lafayette, Indiana       
\endaddress
\email
gottlieb\@math.purdue.edu      
\endemail

\date 
{\thismonth}{\space\number\day,} {\number\year}
\enddate


\abstract
Given a parameterized space of square matrices, the associated set of eigenvectors forms
some kind of a structure over the parameter space. When is that structure a vector bundle? When is there a vector field of eigenvectors? We answer those questions in terms of three obstructions,
using a Homotopy Theory approach. We illustrate our obstructions with five examples. One of those examples gives rise to a 4 by 4 matrix representation of the Complex Quaternions.
This representation shows the relationship of the Biquaternions with low dimensional Lie groups
and algebras, Electro-magnetism, and Relativity Theory. The eigenstructure of this representation
is very interesting, and our choice of notation produces important mathematical expressions found  in  those fields and in Quantum Mechanics. In particular, we show that the Doppler shift factor is
analogous to Berry's Phase.
\endabstract

\define\BB{\bold B}
\def\ov{\overline}

\define\BE{\bold E}

\define\bE{\bold E}
\define\bB{\bold B}

\define\bw{\bold w}

\define\bC{\Bbb C}
\define\bH{\Bbb H}
\define\bR{\Bbb R}
\define\bK{\Bbb K}
\define\bZ{\Bbb Z}



\define\ds{\displaystyle}

\define\thismonth{\ifcase\month 
  \or January\or February\or March\or April\or May\or June%
  \or July\or August\or September\or October\or November%
  \or December\fi}
\NoBlackBoxes

\endtopmatter
\parskip=6pt
\document 

\head 1.\ Introduction\endhead

This work was stimulated by the Gibbs Lecture of Sir Michael Berry given at
the 2002 American Math. Soc. meeting in San Diego California. Berry's 
lecture discussed the discription of physical phenomina by means of 
slowly changing eigenvectors of relevant linear operators, usually 
Hamiltonians of Quantum Mechanics. This work was advanced by several
mathematical physicists, such as Barry Simon, under the name of Berry's
Phase. The original papers are \cite{Berry(1984)} and \cite{Simon(1983)}. A multitude of similar
phenomena are found in \cite{Berry(1990)}.

Berry's Phase can be thought of in terms of {\it eigenbundles}, or 
{\it spectral bundles}
as some mathematical physicists call them. These are vector 
bundles whose fibres are spaces of eigenvectors associated to linear
operators which are parameterized by the base space.
 
There are two questions involving these spectral bundles. The first is:
When do they exist? The second is: What is a relevant connection to 
put on a spectral bundle which results in physical descriptions?

The first question is topological, the second is more geometrical and
of course physical. We will approach the first question from a 
homotopy theoretical point of view. Spectral bundles are related to
an area of Analysis concerned with spectral projections. Mathematical
physicists have incorporated some homotopy concepts, such as homotopy
groups, in their study of spectral bundles, \cite{Avron, Sadun, Segert, Simon(1989)}. What we do here is study
the existence of spectral bundles by means of a commutative diagram.
This will characterize when spectral bundles exist in terms of three
obstructions, and will organize the many variants under which the
existence problem can be posed.

We illustrate the issues involved by giving a few simple examples and
one sophisticated example. The sophisticated example consists of a set
of $4 \times 4$ matrices which are a representation of the {\it biquaternions},
that is the quaternions complexified. We denote the quaternions by 
$\bH$ and the biquaternions by $\bH \otimes \bC$.

The quaternions and biquaternions have been studied for over 150 years
as a convenient language for physics, \cite{Gsponer, Hurni(2002)} 
The generalization of quaternions,
called Clifford Algebras, has also been extensively studied by physicists,
especially by Dave Hestenes under the name of {\it Geometric Algebra}, \cite{Hestenes, Sobcyk(1987)}. 

Our particular 4--dimensional representation of the biquaternions
naturally gives rise to 4--dimensional representations of important low
dimensional Lie groups and algebras. There is a conjugate representation
also, and a ``modulus square mapping'', $\frak m$ , from these representations of the
biquaternions gives  well known relationships of low dimensional
Lie groups, and electromagnetic energy-momentum tensors,
as well as a cononical form of the eigenvectors of Lorentz transformations.
This last feature allows us to see the Doppler shift factor as an analogue of
Berry's phase. Finally in Section 8, we give two examples of probability distributions in
Quantum Mechanics which can be expressed as inner products of eigenvectors .

\head 2.\ Examples \endhead

In this section we set up our basic point of view and illustrate with 4 examples.

Let $V$ be a vector space over the Real numbers  $\bR$  or the Complex numbers
$\bC$. Consider the space $Hom(V,V)$ of linear maps from $V$ to $V$.
We assume that $V$ is a finite dimensional space so that we can describe the
topology of $Hom(V,V)$ simply. If a basis is chosen for the n--dimensional
space $V$, then we have automatically chosen an isomorphism from $Hom(V,V)$ to
$M_n(\bK)$ where $\bK$ stands for either the scalars $\bR$ or $\bC$. 
Here $M_n(\bK)$ denotes
the space of $n \times n$ matrices with entries in $\bK$.
This space is given the Euclidean topology of $\bK^{n^2}$.

Now let $\Phi :B \rightarrow  Hom(V,V)$ be a continuous map where 
$B$ is a topological space. We will call $\Phi$ a {\it field} of
linear operators (or matrices) on $B$. In the physics literature, 
this is frequently called a system of linear operators parametrized
by $B$. In Physics, $B$ is usually an interval of the Real line and
the parameter is frequently thought of as time. Another variant is
the field is over a parameter space $B$, and a physical process is
represented by a path in the parameter space $B$.

There is a trivial example where $\Phi :B \rightarrow  M_n(\bK)$
is the constant map which maps every point to the identity matrix
$I$. In this case, any subbundle of the trivial bundle is an eigenbundle.

At the opposite extreme we give an example for which no eigenbundle exists.
Let $B$ be the rotation group in two dimensions, $SO(2)$, and let 
$\Phi$ be the inclusion 
map of $SO(2)$ into the space of $2 \times 2$ matrices $M_2(\bR)$.
Every rotation except for the identity has imaginary eigenvectors and
eigenvalues, hence there cannot be a real spectral bundle over $SO(2)$.

We will give four examples below which illustrate various issues which
arise in the study of the existence of spectral bundles.

{\bf Example} 1 : Let $B = \bR$ and let $\Phi : \bR \rightarrow M_2(\bR)$ be
given by 

$$ \Phi (t) = \left(\matrix
             1& t\\
             0& 1\\
        \endmatrix \right)
$$

Every $\Phi (t)$ has only one eigenvalue $\lambda = 1$ with corresponding
eigenspace spanned by the vector $(1 , 0)^T$ when $t$ does not equal $0$,
and at $t = 0$, $ \Phi(0) = I$ so the eigenspace is all of $\bR^2$.
In this case the spectral line bundle exists and is trivial since 
there is a nonzero cross-section. For example, the map which takes 
$ t \mapsto (t, (1, 0)^T)$ is a cross-section. We regard this cross-section
as a vector field of eigenvectors.

{\bf Example} 2: Let $B = \bR$ and let $\Phi : \bR \rightarrow M_2(\bR)$ be
given by

$$ \Phi (t) = \left(\matrix
             1&f(t) \\
          g(t)& 1\\
        \endmatrix \right)
$$
where $f(t)$ is a continuous real valued function which 
is greater than zero if $t$ is positive and equal to zero
if $t$ is nonpositive; and $g(t)$ has the opposite property, for example
$g(t) = f(-t)$. In this example again, there is only one eigenvalue 
$\lambda = 1$, but now the eigenspaces are spanned by $(1, 0)^T$ for 
$t > 0$ and $(0, 1)^T$ for $t < 0$, and at $t = 0$ the eigenspace is
$\bR^2$. Thus there is no continuous choice of eigenvectors over $\bR$
and so there is no eigenbundle. However, if we were willing to change the
field $\Phi$ slightly, by letting $f(t)$ be zero in a small interval about
$0$, then we can connect up the $(1, 0)^T$ vector field continously with
the $(0, 1)^T$ vector fields through eigenvectors in $\bR^2$ near $0$.
So example 2 shows that degenerate eigenspaces are an obstruction to
eigenbundles, but under some circumstances, a slight change in $\Phi$
can eliminate the obstruction.

{\bf Example} 3: Let $B = \bR^3$ and let $\Phi : \bR^3 \rightarrow M_2(\bR)$ be
given by

$$ \Phi (u,v,w) = \left(\matrix
             u&v \\
             v& w\\
        \endmatrix \right)
$$

Then $\Phi (u,0,u)$ has only one eigenvalue $\lambda = u$
and the associated eigenspace is the whole of $\bR^2$.
Off the line $l$ given by $\{(u,0,u)\}$ however, $\Phi$ has two distinct
real eigenvalues and the corresponding eigenspaces are one-- dimensional 
and orthogonal, because $\Phi (b)$ is a symmetric matrix.
Let $B' = \bR^3 - l$. Then there are two spectral line-bundles over $B'$.
But neither of them is a trivial line bundle. So there is no eigenvector
field over $B'$.

This is seen by moving around a loop which links the line $l$. The line
bundle over the loop is not trivial, so it looks like a Mobius band. If 
we regard the map $\Phi$ as mapping into $M_2(\bC)$, the eigenbundles 
over $B'$ are complex line bundles and must be trivial since 
complex line bundles are classified by the first Chern class which 
lives in the second cohomology group with integer coefficients. Since
$B'$ is homotopically equivalent to the circle, the second cohomology
group, and hence the Chern class, and hence the line bundle, must be 
trivial.

This example was mentioned by M. V. Berry in \cite{Berry(1990)} on Page 38,
where he states that this phenomenon didn't seem to be widely known in
matrix theory.

The fourth example is more complex, and it is related to the quaternions
$\bH$, the biquaternions $\bH \otimes \bC$, $SL(2, \bC)$, 
$SO(3,1)$, $\frak s \frak o (3,1)$, $SU(2)$ and 
$\frak s \frak u (2)$ and other topics.

\noindent
{\bf Example} 4: Let $B = \bC^3$ and $\Phi:\bC^3 \rightarrow M_4(\bC)$ so that 
$\Phi ( A_1, A_2, A_3)$ is a matrix $F$ such that 

$$
F = \pmatrix
\vbox{
\offinterlineskip\tabskip=0pt
\halign{
\strut # &
# \hfil &
# \vrule \ \ &
\hfil # \hfil &
\hfil # \hfil &
\hfil # \hfil\cr
& 0 & & $A_1$ & $A_2$ & $A_3$\\
\omit&\multispan{5}{\hrulefill}\cr
& $A_1$ & & 0 & $-iA_3$ & $iA_2$\cr
& $A_2$ & & $iA_3$ & 0 & $-iA_1$\cr
& $A_3$ & & $-iA_2$& $iA_1$ & 0\cr
}}\endpmatrix
$$

Or in block form,

$$
 F  = \left(\matrix 0& \vec A^T\\
               \vec A& \times (-i \vec A)\\
            \endmatrix \right)
$$

\noindent 
where the notation $\times (-i \vec A)$ symbolizes the $3 \times 3$ matrix
which operates on a column vector $v$ to produce the cross product
$ v \times (-i \vec A)$.

Let $\cdot$ represent the usual Euclidean inner product extended linearly to
the complex case. Thus $\vec A \cdot \vec A = A_1A_1 + A_2A_2 + A_3A_3$.
Then the eigenspace structure of $\Phi (\vec A)$ depends on 
$\vec A \cdot \vec A$.

{\bf Case} 1: $\vec A \cdot \vec A \neq 0$. In this case there are two 
nonzero eigenvalues, one the negative of the other (since the square
of the eigenvalue equals $\vec A \cdot \vec A$). Each eigenvalue corresponds
to a two--dimensional eigenspace. Let $B_1$ denote the set of all vectors
$\vec A$ such that $\vec A \cdot \vec A \neq 0$. Then there are no 
eigenbundles for $\Phi$ restricted  to $B_1$.

{\bf Case} 2: $\vec A \cdot \vec A =  0$ and $\vec A \neq 0$. In this case there
is only one eigenvalue, $0$, and it corresponds to a two--dimensional 
eigenspace. Let $B_2$ denote the set of all vectors
$\vec A$ such that $\vec A \cdot \vec A =  0$ and $\vec A \neq 0$. 
Then there is an eigenbundle of rank two over $B_2$. It splits as a
Whitney sum of two trivial line bundles. So there are two linearly independent
eigenvector fields over $B_2$, and one of them consists of real eigenvectors.

{\bf Case} 3: $\vec A = 0$. In this case $\Phi (\vec 0)$ is the zero matrix, so
every vector in $\bC^4$ is an eigenvector.

The above assertions are proved in \cite{Gottlieb(1998), (2001)}. See section 7 of this paper.

\head 3.\ Obstructions to the existence of eigenbundles \endhead

We will show that the obstruction to the existence of spectral bundles
over $B$ for the field $\Phi : B \rightarrow Hom(V,V)$ consists of two
crossections which must be constructed over $B$. A {\it cross-section} to a 
continuous map $f : X \rightarrow Y$ is a map $s : Y \rightarrow X$
so that the composition $f \circ s $ is the identity map, $1_Y$, on $Y$. 
This means that we are able to choose in a continuous way one element in
each fibre $f^{-1}(y)$ of $f$. A cross-section is a homeomorphism of $Y$ 
to its image $s(Y)$ in $X$. Thus we may regard $Y$ as a subspace $s(Y)$
of $X$.

If the first two cross-sections, $s_1$ and $s_2$ exist, then the existence
of a third, $s_3$, gives an eigenvector field.

Suppose we want to construct a spectral bundle whose fibres are $k$--dimensional
eigenspaces over a field $\Phi : B \rightarrow Hom(V,V)$ where $V$ is an
$n$ dimensional vector space. Then we first consider the product space 
$B \times \bK  \times  G_{k,n} \times V $. Here $G_{k,n} = G(V)$ is 
the Grassmannian space of $k$--planes in $V$. 

We define a subspace
$L_3$ of $B \times \bK  \times  G_{k,n} \times V $ as follows:
$L_3$ consists of all the points $(b, \lambda, W, \vec v)$ in 
$B \times \bK  \times  G_{k,n} \times V $ so that $\lambda$ is an
eigenvalue of $\Phi (b)$, and $W$ is a $k$--dimensional eigenspace associated
to $\lambda$, and $\vec v$ is an eigenvector in $W$.

Now the projections 
$$
B \times \bK  \times  G_{k,n} \times V \overset \pi_3\to \longrightarrow B \times \bK  \times  G_{k,n} \overset \pi_2\to \longrightarrow B \times \bK \overset 
\pi_1\to \longrightarrow B
$$
give rise to a sequence of mappings
$$
L_3 \overset \pi_3\to \longrightarrow L_2 \overset \pi_2\to \longrightarrow  
 L_1 \overset \pi_1\to \longrightarrow B
$$
where $L_2 := \pi_3(L_3)$ and $L_1 := \pi_2(L_2)$ are the images of the
projections $\pi_3$ and $\pi_2$ respectively. That is: $L_2$ and 
$L_1$ are the subpaces of  
$B \times \bK  \times  G_{k,n}$ 
and $B \times \bK  $ consisting of the points 
$(b,\lambda,W)$ and $(b,\lambda)$ respectively where
$\lambda$ is an
eigenvalue of $\Phi (b)$, and $W$ is a $k$--dimensional eigenspace associated
to $\lambda$.

Now the map $ \pi_3 : L_3 \rightarrow L_2 $ is a $k$--plane vector bundle.
In fact it is a $k$--spectral bundle with respect to the matrix field 
 $ L_2 \rightarrow M_n(\bK)$ defined by 
$(b,\lambda,W) \mapsto \Phi (b)$.
Now this spectral bundle restricts to a subspace as a spectral bundle over
the matrix field restricted to the subspace. So if $s: B \rightarrow L_2$
is a cross-section to the map $ \pi_1 \circ \pi_2 :L_2 \rightarrow B$, then
the restriction of the spectral bundle over $L_2$ to the spectral bundle
over $s(B)$ gives a spectral bundle $\pi_3 : L_3' \rightarrow s(B)$ 
over $B$ for the matrix field $\Phi$. 
 
The above paragraphs give the notation and the proof for the following
classification theorem for spectral bundles:

\proclaim{Theorem 3.1}The $k$-spectral bundles are in one to one 
correspondence with the cross-sections of the map 
$ \pi_1 \circ \pi_2 :L_2 \rightarrow B$
\endproclaim

It is convenient to break the cross-section $s$ into two cross-sections:
$s_1 : B \rightarrow L_1$, and $s_2: s_1(B) \rightarrow L_2'$ where $L_2'$
denotes $\pi_2^{-1}(s_1(B))$, the preimage of $s_1(B)$ contained in $L_2$.
Now the composition $s_2 \circ s_1$ is a cross-section to 
$ \pi_1 \circ \pi_2 :L_2 \rightarrow B$. On the other hand, a cross-section
$s : B \rightarrow L_2$ induces the cross-section $\pi_2 \circ s =: s_1$,
and the cross-section $s_2$ is
$s \circ \pi_1 : s_1(B) \rightarrow L_2'$ .

The following diagram may be helpful in tracing the above notation in the theorem
below. The horizontal arrows represent inclusion maps.

$$
\CD
L_3'' @>>> L_3'@>>> L_3 @>>> B \times \bK  \times  G_{k,n} \times V\\  
@As_3AA    @VV\pi_3V    @VV\pi_3V @V\pi_3VV \\
s_2s_1B @>>>  L_2' @>>> L_2@>>> B \times \bK  \times  G_{k,n} \\
@.    @As_2AA    @VV\pi_2V@V\pi_2VV \\
@.s_1 B @>>> L_1 @>>>B \times \bK  \\
@.    @.   @As_1AA @V\pi_1VV \\
 @.  @. B @=B
\endCD 
$$

\proclaim{Theorem 3.2} 
\item{a)} The set of $s_1$ cross-sections is in one to one
correspondence with the continuous functions $ \lambda :B \rightarrow \bK$ 
so that every every $\lambda (b)$ is an eigenvalue of $\Phi (b)$ whose 
associated eigenspace has dimension $\geq k$.
\item{b)} The set of $s_2$ cross-sections 
corresponds to the continuous selections
of $k$--dimensional subspaces of eigenvectors with eigenvalues $\lambda (b)$.
\item{c)} The set of nowhere zero cross-sections $s_3$ of the spectral
bundle $L_3'' \overset \pi_3\to \longrightarrow s_2s_1(B) = B $ corresponds
to the set of nowhere zero eigenvector fields for the eigenbundle.
\endproclaim

\demo{Proof}
\item{a)}The cross-section $s_1(b) = (b, \lambda (b))$ is continuous if and only if
$\lambda (b)$ is continuous.
\item{b)} $s_2 (b) = (b, \lambda (b), W_b)$ where $b \mapsto W_b$ picks out
a  $k$--dimensional subspace of eigenvectors with eigenvalue $\lambda (b)$
contained in $V$, that is it is a function from $B \rightarrow G(V_k)$. Now
$s_2 $ is continuous if and only if the function $B \rightarrow G_k(V)$ 
is continuous.
\item{c)} $s_3$ is a cross-section to the vector bundle 
$L_3'' \overset \pi_3\to \longrightarrow s_2s_1(B) = B $, so $s_3(b)$ is
an eigenvector for $\Phi (b)$. If $s_3(b) \neq 0$ for all $b$ in $B$, then
the spectral bundle has a trivial line bundle summand, or equivalently, a nonzero
eigenvector field.
\qed
\enddemo

Now let us consider $L_1$ for complex spectral line bundles. This is the
largest of the possible $L_1$'s for a fixed $\Phi$. Every other $L_1$ for
higher dimensional complex spectral bundles, or for real spectral bundles
associated to $\Phi$, must be a subspace of the $L_1$ for complex 
spectral line bundles. In those cases it is possible that there are no
eigenvalues for $\Phi (b)$ and hence there is no cross-section $s_1$.
Examples like the real rotation matrices $SO(2)$ or the spectral 3-bundles
of example 4 show that there is no $s_1$ because $\pi_1$ is not onto.
But for complex spectral line bundles, not only must $\pi_1$ be onto,
but $L_1$ is a topological branched covering of $B$, where we mean the following by 
{\it topological branched covering}:
A space $X$ which admits a continuous onto map $p : X \rightarrow B$ such that all fibres 
are discrete and so that the path lifting property holds. That is for every $x \in X$,
and path $\sigma$ in $B$ starting at $\sigma (0) = p(x)$, there is a path $\overline \sigma$
in $X$ so that $\sigma = p \circ \overline \sigma$ and $\overline \sigma (0) = x$.

\proclaim{Theorem 3.3} For complex line bundles,
$ \pi_1 :L_1 \rightarrow B$
is a topological branched covering of $B$.
\endproclaim

\demo{Proof}
Consider the mapping from $B$ to the complex polynomials of degree $n$
given by $b \mapsto det(\lambda I - \Phi (b))$
This is a continuous map from $b$ to the characteristic polynomial
of $\Phi (b)$. The Fundamental Theorem of Algebra tells us that there are
$n$ roots of this polynomial counting multiplicities, for any point $b$.
The roots are of course, the eigenvalues of $\Phi (b)$. I like to think of
it using vector fields. Over each $b$ in $B \times \bC$ is a fibre $\bC$.
On each fibre there is a vertical vector field on $\bC$ given by attaching
the vector $p_b(z)$ to $z$ where $p_b$ is the characteristic polynomial
for $\Phi (b)$. Each zero has a positive vector field index, equal to the
multiplicity of the corresponding root. The sum of the local indices adds up
to a global index $n$ for every fibre. The set of the the zeros is $L_1$.
So every $b$ is covered by at least one zero and at most $n$ zeros.
Hence $\pi_1$ is onto, and $L_1$ consists of at most $n$ connected components
over $B$. As we move from one $b$ to a nearby point, there are zeros in 
the new fibre close to where they were at $b$, because no zero can be
annihilated by another since there are no nonpositive indices to cancil out.
This gives $L_1$  the branched covering structure. See \cite{Gottlieb, Samaranayake(1994)} for a detailed discussion of
the index of vector fields.
\qed
\enddemo

In the case of real matrices, the real characteristic polynomial
$det(\lambda I - \Phi (b))$ can be thought of as a vertical vector field
on the fibres $\bR$. Again the zeros of this vertical vector field on
$B \times \bR$ gives us $L_1$, but here it is not necessarily a branched
cover over $B$. The reason is that the zeros of the characteristic 
polynomial on the real line have indicial values of $1$, $-1$ or $0$ .
The opposite signs and zero indices allow the zeros on the Real line
to annihilate each other, so that there may not be a nearby zero on a 
nearby fibre to continue the local covering of $B$ by $L_1$.

The total index on each fibre $\bR$ is $1$ for odd order matrices and
$0$ for even order matrices, so the sum of the local indices of each zero
add up to $0$ in even dimensions and $1$ in odd dimensions. Thus, for
odd dimensional matrix fields, there is always a zero of index $1$ in 
each fibre, so $\pi_1$ is always onto in that case. For the even dimensional
matrix field however, there is no guarantee of a zero in every fibre, so
$\pi_1$ may not be onto.

The real matrix field may be considered as acting on a complex vector space.
In this case, the zeros on the real line in $\bC$ still have their indices of
positive integers as well as their indices $\pm 1$ or $0$ on the Real
line. In this case, a real zero's annihilation actually is given by a 
splitting of the zero into two complex conjugate zeros, which of course
are off the Real line. Thus a real zero doesn't disappear, it splits into
two conjugate zeros which leave the Real line in the Complex plane.

Now we will reconsider our examples in light of the above considerations.

{\bf Example} 1 has only one eigenvalue for each $b \in \bR$, so $s_1$ exists.
At each point $b$ there is only one $1$--dimensional eigenspace except at
$b = 1$, where it is $2$--dimensional.
This potentially blocks the existence of $s_2$, but it happens that we may
choose a $1$--dimensional eigenspace in the $2$--dimensional eigenspace so that
the choice of $1$--dimensional eigenspaces is continuous. 
So $s_2$ exists. There is an obvious eigenvector field, so $s_3$ exists. It
is worth remarking that given a vector bundle over a contractible space such
as $\bR$, the vector bundle must be trivial and there are always nonzero
vector fields; or to say it another way, we can always split off a trivial
line bundle. 

{\bf Example} 2 is the same as Example 1, except that it is impossible to choose
a 1--dimensional subspace at $b = 0$ in such a way to make a continuous
selection of 1--dimensional eigenbundles. Hence $s_2$ does not exist. 
The possibility was mentioned of altering $\Phi$ slightly to eliminate
this obstruction to $s_2$ existing.
For 1--dimensional $B$'s such as a line interval or a circle, this can always
be done.
Of course, since $Hom(V,V)$ is contractible,we can always homotopy 
$\Phi$ to a constant and obtain a new $s_2$, but this is too large a 
change for most purposes. 

There are homotopy obstructions to changing 
$\Phi$ so as to eliminate 
the obstruction to $s_2$. Suppose that $D$ is the unit disk in the plane.
Let $B = D$, and let $\Phi (b)$ be a symmetric matrix of order $2$
with eigenvalues $\pm 1$ when $b \in S^1$, where $S^1$ is the boundary of
$D$. Suppose that the $+1$ eigenvectors are pointing orthogonally
outside of $D$. The it is impossible to extend $\Phi$ over $D$ with values
symmetric matrices such that every matrix has no $2$--dimensional eigenspace.
This follows since the outward pointing eigenvector field cannot be extended
to a nonzero vector field over $D$, since such a vector field has
$index = 1$. Since every symmetric matrix has a two frame of eigenvectors
whenever the two eigenvalues are distinct, such an extension 
of $\Phi$ would give rise to a a nonzero vector field. Contradiction. 

{\bf Example} 3 exhibits some homotopy type features. Recall
$$ \Phi (u,v,w) = \left(\matrix
             u&v \\
             v& w\\
        \endmatrix \right)
$$
Since the matrices are
symmetric, the eigenvalues are real and we can find continuous eigenvalue
functions on $B = \bR^3$. Hence $s_1$'s exist. On the other hand, 
$s_2$ does not exist. We know that if an $s_2$ existed, there would be a
eigenbundle over $\bR^3$, which is contractible. Hence it would be a trivial
line bundle. But we know that on a circle linking $l$, the restriction
line bundle is not trivial. So that contradicts the triviality of a bundle
over $\bR^3$. If we consider the question over $B' = \bR^3 - l$, we have
eliminated degenerate eigenspaces, every eigenspace is $1$--dimensional,
so we can choose a continuous selection of eigenspaces, so $s_2$ exists,
and we have a spectral line bundle over $B'$. But it is not a trivial 
bundle. Now real line bundles are classified by their Stiefel-Whitney 
class $w_1$, which lives in the first cohomology group of $B'$ with
$\bZ_2$ coefficients, $H^1(B', \bZ_2)$.
Now $B'$ is homotopy equivalent to $S^1$, and so there is only one nonzero
$w_1 \in H^1(B', \bZ_2) = \bZ_2$.

If we consider the same field acting on a complex two--dimensional vector
space, we again get a spectral line bundle over $B'$, but this time the
bundle is trivial in that is there is a nonzero eigenvector field, but it is
not completely real. A complex line bundle is classified by its Chern class
$c_1 \in H^2(B', \bZ)$, the two--dimensional cohomology group with integer
coefficients. Since $B'$ is homotopy equivalent to a circle, the 
two--dimensional cohomology
must be zero and hence $c_1 = 0$, so the bundle is trivial.

{\bf Example} 4 has the property that every eigenspace has complex dimension $2$
except for the $0$ matrix. If we remove the $0$ matrix from consideration,
we see that if $s_1$ exists, then $s_2$ would exist and we would have an
eigen $2$-bundle. If we restrict to Case 2, the set $B_2$ of all vectors
$\vec A$ such that $\vec A \cdot \vec A =  0$ and $\vec A \neq 0$,
we get $s_1$ since the only eigenvalue is $0$. Hence in this case there
exists an eigenbundle of rank $2$ over $B_2$. Let us write
$\vec A := \vec E + i \vec B$ where $\vec E $ and  $ \vec B$ are real
vectors. In this case, where $\vec A \cdot \vec A = 0$, we have
$E = B$ and $\vec E \cdot \vec B = 0$. We may describe the eigenspace
by means of two linearly independent eigenvectors:
$\vec E + i \vec B$ and
$E^2u + \vec E \times \vec B$ where $u = (1, 0, 0, 0)$. Here we are 
regarding the $3$-vectors as living in the space orthogonal to $u$.
These eigenvectors each give rise to an eigenvector field which shows
that over $B_2$ the eigenbundle of rank $2$ splits as a Whitney sum
of two trivial spectral line bundles.

In Case 1 of Example 4, where $B_1 $ is the set of vectors $\vec A$
such that $\vec A \cdot \vec A \neq 0$, we see that $s_1$ does not exist.
In this case each matrix has two distinct eigenvector spaces.
Recall that for complex line bundles, Theorem 3.3 states that
$ \pi_1 :L_1 \rightarrow B$
is a branched covering of $B$. If we restrict ourselves to matrices so
that every eigenvalue is distinct, then the branching part of the branched
covering is eliminated and we have a covering. Each connected component
of the covering space is a connected covering space. A cross-section
$s_1$ exists if and only if there is a connected component which is 
homeomorphic to $B$, that is, if and only if there exists a one to one
covering of $B$. In situation at hand, the eigenvalues are 
are not distinct, but there are only two of them, one being the 
negative of the other. This gives rise to a two to one 
covering of $B_1$. Hence $s_1$ does not exist.

In this case, if we move around a closed curve in $B_1$ which loops
$B_2$ one time, we arrive at the same matrix, but the eigenspace
has been transported to the eigenspace corresponding to the opposite
eigenvalue. This is a subtle effect when encountered without the aid
of the double covering point of view.

We will add one more example to our list of four examples. This will
actually be an extension of Example 4, and is a faithful $4$--dimensional
representation of the Biquaternions $\bH \otimes \bC$.

{\bf Example} 5: Consider the set $I+S$ of all $4 \times 4$ matrices of the form
$aI + F$ where $a$ is any complex number and $I$ is the identity matrix
and $F$ is any matrix from Example 4. That is $F \in S$ and so has the form
$$
F  = \left(\matrix 0& \vec A^T\\
               \vec A& \times (-i \vec A)\\
            \endmatrix \right)
$$

Here $B = \bC^4$, and $\Phi (A_0,A_1,A_2,A_3) = A_0I + F$. That is:
$$
\Phi (A_0,A_1,A_2,A_3) = \pmatrix
\vbox{
\offinterlineskip\tabskip=0pt
\halign{
\strut # &
# \hfil &
# \vrule \ \ &
\hfil # \hfil &
\hfil # \hfil &
\hfil # \hfil\cr
& $A_0$ & & $A_1$ & $A_2$ & $A_3$\\
\omit&\multispan{5}{\hrulefill}\cr
& $A_1$ & & $A_0$ & $-iA_3$ & $iA_2$\cr
& $A_2$ & & $iA_3$ & $A_0$ & $-iA_1$\cr
& $A_3$ & & $-iA_2$& $iA_1$ & $A_0$\cr
}}\endpmatrix
$$

Let $\langle \  , \rangle $ represent the usual Minkowskian inner product extended linearly to
the complex case. Thus, if $A :=  (A_0,A_1,A_2,A_3) =: (A_0, \vec A)$, then
$\langle  A ,  A \rangle= -A_0A_0 + A_1A_1 + A_2A_2 + A_3A_3 = -A_0A_0 + \vec A \cdot \vec A$.
Then the eigenspace structure of $\Phi ( A)$ depends on
$\langle  A ,  A \rangle$.

{\bf Case} 1: $\langle  A ,  A \rangle \neq 0$ and $\vec A \neq 0$. In this case there are two
nonzero eigenvalues when $\vec A \cdot \vec A \neq 0$.
Each eigenvalue corresponds
to a two--dimensional eigenspace. Let $B_1$ denote the set of all vectors
$ A$ such that $\langle  A ,  A \rangle \neq 0$. Then there are no
eigenbundles for $\Phi$ restricted  to $B_1$.

{\bf Case} 2: $\langle  A ,  A \rangle =  0$ and $\vec A \neq 0$. 
In this case there is one or two eigenvalues, but one of them is equal
to $0$, 
and it corresponds to a two--dimensional
eigenspace. Let $B_2$ denote the set of all vectors
$ A$ such that $\langle  A ,  A \rangle =  0$ and $\vec A \neq 0$.
Then there is an eigenbundle of rank two over $B_2$. It splits as a
Whitney sum of two trivial line bundles. So there are two linearly independent
eigenvector fields over $B_2$, and one of them consists of real eigenvectors.

{\bf Case} 3: $ \vec A = 0$. In this case $\Phi ( A)$ is a diagonal matrix, so
every vector in $\bC^4$ is an eigenvector.

We note that the cases of Example 5 seems to be very similar to the cases of  Example 4, 
but now the 
eigenvalues are not each other's negatives, and in Case 2 there are one
{\it or} two eigenvalues. But one of them is always zero, so $s_1$ exists in that
case since the eigenvalue map is the constant zero. But then the nonzero
eigenvalue also must form an eigenfunction over $B_2$, and so there is another
spectral $2$-bundle over $B_2$. Over the region where $A_0 = 0$, this
second spectral $2$-bundle is identical with the first.

\head 4.\ Biquaternions \endhead

The set of matrices of Example 5 
$$
\Phi (A_0,A_1,A_2,A_3) = \pmatrix
\vbox{
\offinterlineskip\tabskip=0pt
\halign{
\strut # &
# \hfil &
# \vrule \ \ &
\hfil # \hfil &
\hfil # \hfil &
\hfil # \hfil\cr
& $A_0$ & & $A_1$ & $A_2$ & $A_3$\\
\omit&\multispan{5}{\hrulefill}\cr
& $A_1$ & & $A_0$ & $-iA_3$ & $iA_2$\cr
& $A_2$ & & $iA_3$ & $A_0$ & $-iA_1$\cr
& $A_3$ & & $-iA_2$& $iA_1$ & $A_0$\cr
}}\endpmatrix
$$
is a representation of the biquaternions.

Obviously it is isomorphic to $\bC^4$ as a vector space. We will list a 
basis below which will reveal the relationship of the matrices and the
biquaternions. Let $x$ denote the matrix above in which $A_1=1 \text { and
the other } A_i = 0$.

That is 
$$
x = \Phi (0,1,0,0) = \pmatrix
\vbox{
\offinterlineskip\tabskip=0pt
\halign{
\strut # &
# \hfil &
# \vrule \ \ &
\hfil # \hfil &
\hfil # \hfil &
\hfil # \hfil\cr
& $0$ & & $1$ & $0$ & $0$\\
\omit&\multispan{5}{\hrulefill}\cr
& $1$ & & $0$ & $0$ & $0$\cr
& $0$ & & $0$ & $0$ & $-i$\cr
& $0$ & & $0$& $i$ & $0$\cr
}}\endpmatrix
$$

In the same way we define matrices 

\noindent $y := \Phi (0, 0, 1, 0)$

\noindent $z:= \Phi (0, 0, 0, 1)$

\noindent $I = \Phi (1, 0, 0, 0)$, the identity matrix of order $4$.

Now $xy = iz$ and $ x^2 = y^2 = z^2 = I$ and $ xy = -yx$.
Then the basis $\{ ix,\ \ iy,\ \ iz,\ \ I \}$ 
obviously has the relations defining the biquaternions.

There is another representation of the biquaternions in which the traceless matrices
are given by 
$$
F  = \left(\matrix 0& \vec A^T\\
               \vec A& \times (i \vec A)\\
            \endmatrix \right)
$$
These matrices differ from the previous set in Example 4 by changing
the $-i$ to $+i$. If we denote the set of matrices of Example 4 by $S$,
let $ \overline S$ denote the set of matrices of the form $F$ above.

Now let $\{ X,\ \ Y,\ \ Z,\ \ I \}$ 
be the complex conjugates of $\{ x,\ \ y,\ \ z,\ \ I \}$ respectively.
These new elements satisfy
 $XY = -iZ$ and $ X^2 = Y^2 = Z^2 = I$ and $ XY = -YX$.
So the basis $\{ -iX,\ \ -iY,\ \ -iZ,\ \ I \}$
obviously has the relations defining the biquaternions for 
$I \oplus \overline S$. 

Now it happens that any $F \in S$ commutes
with any $G \in \overline S$. That is $F \overline G = \overline G F$ for
$F, G \in S$. This gives rise to a pairing
$(I \oplus S) \otimes (I \oplus \overline S) \rightarrow M_4(\bC)$
given by $A \otimes B \mapsto AB$ where the product $AB$ is in the space of
$4 \times 4$ complex matrices.
This pairing is an isomorphism of rings. This can be seen by observing that
the following set of sixteen matrices forms a basis of $M_4(\bC)$:
\proclaim{Theorem 4.1} The set of sixteen matrices
$$
\matrix\format\r&\quad\r&\quad\r&\quad\r&\quad\r\\
I,   & x X,& y Y,& zZ,\\
x,&X,&yZ,&zY,\\
y,&Y,&xZ,&zX,\\
z,&Z,&xY,&yX
\endmatrix\tag"a)"
$$
forms a basis for $M_4(\bC)$, the vector space of $4\times 4$ complex matrices.

b) The square of each of the matrices in the basis is $I$.

c) Each matrix is Hermitian, so real linear combinations of the basis are the
$4 \times 4$ Hermitian matrices.

d) Every matrix has zero trace except for $I$.
\endproclaim
\demo{Proof}
Theorem 3.3 of \cite{Gottlieb(2001)}.
\enddemo
It is easy to calculate any $4 \times 4$ matrix in terms of this basis
using MATLAB.
 Below I produce a matrix whose first column is  x  written as a column
 vector of length 16. (this is done by x(:), which counts from 1 down
 the first column and then down the next column until you arrive at
 the $4 \times 4$ term which is the last number of the vector).
  The remaining columns are given in the order as shown below in the
 definition of  Total.
 
\noindent (4.1)

\noindent Total = [x(:) X(:) y(:) Y(:) z(:) Z(:) xY(:) yX(:) yZ(:) zY(:) zX(:) xZ(:) xX(:) yY(:) zZ(:) I(:)];

Now any $4 \times 4$ matrix M can be converted into a vector M(:).
The command  Total $\backslash$ M(:) gives the vector of coefficients which when
multiplied with the basis in the order found in Total above will give the
linear combination of M in terms of the basis.

Now since $M_4(\bC)$ is the complex Clifford algebra $\Cal C \ell(4)$, 
there must be generators $\alpha_0,\alpha_1,\alpha_2,\alpha_3$ so that 
$\alpha_i\alpha_j+ \alpha_j\alpha_i=\delta_{ij} I$.
One such set of $\alpha$'s is given by
$$
\align
\alpha_0&=x\\
\alpha_1&=y\\
\alpha_2&=zX\\
\alpha_3&=zY.
\endalign
$$

\proclaim{Theorem 4.2} 
Let $F, G \in S$ satisfy $Fu = \vec A$ and $Gu = \vec B$. Then

a) $FG + GF = (\vec A \cdot \vec B)I$

b) $F \overline G = \overline G F$

c) $[F,G]u = 2i\vec A \times \vec B$

d) $e^{F}=\cosh(\lambda_{F}) I+{\sinh(\lambda_{F})\over \lambda_{F}}\ F$
where $\lambda_{F}$ is an eigenvalue of $F$.
\endproclaim
\demo{Proof}
Corollary 4.7, Theorem 4.8, Corollary 4.4, and Theorem 8.5 of \cite{Gottlieb(1998)} respectively.
\qed
\enddemo

Now every  nonsingular 
matrix $A \in M_4(\bC)$ gives rise to an inner automorphism of
$M_4(\bC)$ given by $B \mapsto A^{-1}BA$.
These maps transform the basis into a new basis with the same algebraic
properies, but the form of the representative matrices can be quite
different. We will end this section discussing what distinguishes our representation
from the other representations.

The matrices of $S$ (or $\overline S$) are skew symmetric with respect
to the Minkowski metric $-+++$. That is equivalent to the property
$F^T = - \eta F \eta$ where $F \in S,\  \ F^T$ is the transpose of $F$ and
$\eta = $ the diagonal matrix with $-1, 1, 1, 1$ down the main diagonal.
A popular set of matrices are the skew symmetric matrices with respect to the 
Euclidean metric. They satisfy $A = - A^T$.
Now $\eta^{1/2} F \eta^{-1/2}$ is a skew symmetric matrix if
$\eta^{1/2}$ and $ \eta^{-1/2}$ equal the diagonal matrix $\pm i, 1, 1, 1$
respectively. Hence if
$$ F =\pmatrix
\vbox{
\offinterlineskip\tabskip=0pt
\halign{
\strut # &
# \hfil &
# \vrule \ \ &
# \hfil\cr
& 0 & & $\vec A^T$\\
\omit&\multispan{3}{\hrulefill}\cr
& $\vec A$ & & $\mp i(\times \vec A)$\cr
}}\endpmatrix \tag4.2$$ 
then
$$ \eta^{1/2} F \eta^{-1/2} = -i \pmatrix
\vbox{
\offinterlineskip\tabskip=0pt
\halign{
\strut # &
# \hfil &
# \vrule \ \ &
# \hfil\cr
& 0 & & $- \vec A^T$\\
\omit&\multispan{3}{\hrulefill}\cr
& $\vec A$ & & $\pm (\times \vec A)$\cr
}}\endpmatrix \tag4.3$$ 
Thus $M_4(\bC)$ is the tensor product 
$(I + \eta^{1/2} S \eta^{-1/2}) \otimes 
(I + \eta^{1/2} \overline S  \eta^{-1/2})$.
This means that the transformed $S$ matrices still have squares equal to a multiple
of the identity, and it satisfies the same exponential equation as in Theorem 4.2d. And the
transformed $S$ and $\overline S$ still commute, but they are no longer the complex conjugate
of each other. It is this property which gives our representation its distinctive advantage, because
the "modulus squared map" is a multiplicative homomorphism on $S$.

The matrices of the form
$$ F =\pmatrix
\vbox{
\offinterlineskip\tabskip=0pt 
\halign{
\strut # &
# \hfil &
# \vrule \ \ &
# \hfil\cr
& 0 & & $\vec A^T$\\
\omit&\multispan{3}{\hrulefill}\cr
& $\vec A$ & & $\times \vec C$\cr
}}\endpmatrix \tag4.4$$
are the skew symmetric matrices with respect to the Minkowski inner product.
So $S$ and $\overline S$ 
are skew symmetric matrices with respect to the Minkowski inner product.
The only skew symmetric matrices with respect to the Minkowski inner product
whose squares are multiples of $I$ are precisely the matrices of
$S$ and $\overline S$. \cite{Gottlieb(1998)}, see Theorem 4.5 .

Now note that if $F \in S$, then both the complex conjugate $\overline F$
and the transpose $F^T$ are both in $\overline S$. Thus the pseudo 
automorphisms $conjugation: A \mapsto \overline A$, which is antilinear
in that it changes the sign of $i$, and $transpose: A \mapsto A^T$, which
reverses the order of multiplication, interchange $S$ and $\overline S$. In
terms of our basis,
$ax + by + cz \mapsto \overline aX + \overline  bY + \overline  cZ$
under conjugation and
$ax + by + cz \mapsto aX + bY + cZ$ under transposition. 
The composition of conjugation and
transposition yields the Hermitian conjugate 
$\dagger : ax + by + cz \mapsto \overline ax + \overline  by + \overline  cz$
which is an antilinear isomorphism which preserves $S$ and $\overline S$.

On the other hand, $S$ and $\overline S$ are interchanged 
by the inner automorphism $A \mapsto \eta A \eta$. That follows since
$\eta F \eta = -F^T$ when $F \in S$. In terms of our basis,
$ax + by + cz \mapsto -aX - bY - cZ$. 

\head 5.\ The Modulus squared map \endhead

We define the modulus squared map and list several of its properties in this section.
\proclaim{Definition}The modulus squared map is a multiplicative homomorphism
 ${\frak m}: (I + S) \rightarrow M_4(\bR)$
given by $ A \mapsto \frak m(A) = \overline A A$. Its image $\frak m (I + S)$ is denoted 
by $\frak M$.
\endproclaim

To show that this definition is well--defined, we must show that its image is in the set of real 
matrices; and that it preserves matrix multiplication. The following lemma does that.

\proclaim {Lemma 5.1}
Suppose $A$ and $B$  square matrices. Then 
\item{a)}$\overline AA$ is a real matrix if and only if
$A$ and $\overline A$ commute.
\item{b)} $\frak m (AB) = \frak m (A) \frak m (B)$
\endproclaim
\demo{Proof}
\item{a)} A matrix is real if and only if it is equal to its own complex conjugate.
Now 
$A \overline A = \overline A A = \overline{A \overline A}$ since $A$ and $\overline A$ commute.
Conversely, suppose $A \overline A$ is real. Now $A = C + iD$ 
where $C$ and $D$ are real. So
$A \overline A = (C + iD)(C - iD) = C^2 - D^2 + i[D,C]$.
Since $A \overline A$ is real, the commutator $[D,C] = 0$. This implies
that $A \overline A = \overline A A$.
\item{b)} First of all , note that $I+S$ is closed under multiplication. See Lemma 6.2. Then
 $\frak m (AB) = AB \overline AB = A \overline A B \overline B =  \frak m (A) \frak m (B)$.
\enddemo

We will call $\frak m$ the modulus squared map in analogy with the 
complex absolute value squared of a complex number.

Now $\frak m$ has many striking properties. The following are the most 
interesting.

\proclaim{Theorem 5.2}
The set $\frak M$ is homeomorphic to the cone over the projective space
$PC^3$
\endproclaim
\demo{Proof}
As a vector space $(I + S)$ is isomorphic to $\bC^4$. The modulus map $\frak m$ has
fibres $S^1$ over all points of $\frak M$ (except for $0$) since $\frak m(F) = \frak m(\alpha F)$
when $\alpha$ is a complex number of unit modulus. Then $\frak m$ can easily be seen to
be an identification map, and the identification of $\bC^4$ by identifying any vector to its
multiple by a scalar with the same modulus is the cone over $PC^3$ with $0$ as the vertex
of the cone.
\enddemo

\proclaim{Corollary 5.3}
The image of $\frak m$ restricted to the unit $7$-sphere in $I + S$ is
the complex projective space $PC^3$
\endproclaim

The {\it Lorentz group} is the set of linear transformations $L$ on Minkowski space which preserves
the Minkowski metric, that is $<Lu,Lv>=<u,v>$. It has four connected components. The 
component containing the identity is called the proper Lorentz group and is denoted by
$SO^+(3,1)$. 

The {\it complex Lorentz group} is the set of linear transformations on complexified Minkowski
space $\bR^{3,1} \otimes \bC$ which preserve in Minkowski metric. The
complex Lorentz group, $L(\bC)$, has two connected components. It plays a role in physics,
\cite{Wightman(2000)}.

The identity component of the complex Lorentz group intersects $I + S$ in a subgroup, which
I will call the {\it biquaternion Lorentz group}. Similarly, the identity component of the Lorentz 
group intersects $I + \overline S$ in a subgroup which is isomorphic to the other by compex
conjugation. The other complex Lorentz group component is disjoint from both biquaternions.

\proclaim{Theorem 5.4}
The image of $\frak m$ restricted to the biquaternion Lorentz group,  which consists of the set
$\{ aI + F  \ \vert \ a^2I - \lambda^2 = 1 \}$, is the real proper Lorentz group
$SO^+(3,1)$. 
\endproclaim

\proclaim{Corollary 5.5}
The Lorentz Group $SO^+(3,1)$ is exponential, that is it has a 
surjective exponential map from $\frak s \frak o^+(3,1)$ .
\endproclaim

We will prove Theorem 5.4 and Corollary 5.5 in the next section. Corollary 5.5 was proved
in \cite{Nishikawa (1983)}. In fact Nishikawa shows that $SO(n,1)$ is exponential.

\proclaim{Theorem 5.6}
$\frak m(S)$ = The set of electromagnetic energy-momentum tensors.
\endproclaim
\demo{proof}
Suppose $F \in S$. Then $Fu = \bold E + i \bold B $ , and if we imagined $\bold E$ and
$\bold B$ as electric and magnetic vectors, then the corresponding electro-magnetic tensor
$T = {1 \over 2} F\overline F$. See Proposition 5.1  with Definition 3.8 in \cite{Gottlieb(1998)} . See 
\cite{Parrott (1987)} for a mathematical account of electro-magnetic energy-momentum tensors.
\enddemo

\proclaim{Theorem 5.7}
$\frak m(S^3) = SO(3)$ where $S^3$ is the unit 3-sphere, that is the 
real unit quaternions.
\endproclaim
\demo{proof}
The real unit quaternions are represented by 
$\{aI +bix+ ciy +diz)$ where $x, y, z$ are the basis matrices of section 4, and $a, b, c, d$
satisfy $a^2 + b^2+c^2+d^2=1$ and are real numbers. If we multiply $\{aI +bix+ ciy +diz)$
by a unit modulus complex number, the element remains in the real quaternions if and only
if the number is $\pm 1$. Thus $\frak m$ is a 2--1 covering map, so its image must be
$SO(3)$. \qed
\enddemo

The real unit quaternions $S^3$ acts on the right of unit biquaternions
$S^7 =\{ aI +bx + cy +dz \vert a\overline a + b\overline b + c\overline c + d\overline d  = 1 \}$. 
The quotient map is the famous Hopf fibration $S^3 \rightarrow S^7 \rightarrow S^4$. Now
$\frak m : S^7 \rightarrow CP^3$ is a principal $S^1$-fibre bundle and is an equivariant  
map from the free $S^3$ action on $S^7$ to the induced $SO(3)$ action on $CP^3$.
The action of $SO(3)$ on $CP^3 $ is not free.  

Consider the set of matrices in $1 + S$ of the form 
$\{ aI + F \ | \ a^2 = \lambda^2 \}$, where $\lambda$ is the eigenvalue of $F$.
These matrices are those $aI + F$ such that
$(aI + F)(aI - F) = 0$. In biquaternion jargon, these are called {\it nullquats}
or singular quaternions.
Since $ F(\lambda I + F) = \lambda (\lambda I + F)$, we see that the image of $\lambda I + F$
consists of the eigenvectors of $F$ corresponding to the eigenvalue 
$\lambda $. The fact that 
$(\lambda I + F)(\lambda I - F) = 0$ implies that the kernel of
$F(\lambda I + F)$ consists of the eigenvalues of $F$ corresponding to
$- \lambda $. Thus $\lambda I + F$ has rank two. But it is not a spectral
projection unless $\lambda = 1/2$. When $\lambda = 0$ we have the {\it null matrices}
$N$ such that $N^2 = 0$. Here the eigenvector space is both the image and the
kernel of $N$. So $N$ cannot be made into a projection by scalar multiplication. 
However, $N$ does
map $\bC^4$ onto the subspace of eigenvectors of $N$. 

\proclaim{Theorem 5.8}
The image of a nullquat under $\frak m$ is a linear transformation from $\bR^4$ to a real
null 1-dimensional subspace of eigenvectors of the nullquat.
\endproclaim
\demo{proof}
See Theorem 6.7c in \cite{Gottlieb(1998)}.
\enddemo

\head 6. The Exponential Map
\endhead

In this section we show that the exponential map for the proper Lorentz group is surjective
using novel methods.

In order to discuss  eigenvector spaces and exponential maps more fully, we will change
our notation to emphasize the real matrices. We shall follow the notation of 
 \cite{Gottlieb(1998) and (2001)}. 

Let $F \in S$ now be denoted by $cF$ where
$$
cF:=\pmatrix 0 & \bold  A^T\\ \bold  A & \bold \times (-i\bold  A)\endpmatrix
\ \ \text { where } \ \
\bold  A=\bold E +i\bold  B 
$$

Then $cF := F - iF^*$ where $F$ now denotes the real part of $cF$ and $-F^*$ is the
imaginary part. Thus
$$
F=\pmatrix 0 & \bold E^T\\  \bold E & \times\bold  B\endpmatrix \ \text{ and }\
F^*=\pmatrix 0 & -\bold  B^T\\  -\bold  B & \times\bold E\endpmatrix. 
$$

Similarly we define $\overline c F :=  F +  iF^*$. 

Now $F$ is a linear transformation on $\bR^4$ which is skew symmetric with
respect to the Minkowski metric, and $cF$ will be called its {\it complexification }. We may
regard $F$ as a 1-1 tensor corresponding to a two-form $\hat F$. Then
$F^*$ corresponds to the Hodge dual $*\hat F $. If we apply the modulus 
squared map to $cF$, we get $\overline cF cF := 2T_F$ where $T_F$ has the
form of a multiple of the energy-momentum tensor of the electromagnetic
field two-form $\hat F $ corresponding to $F$. On the other hand we may regard
$F$ as an element of the Lie algebra $\frak s \frak o (3,1)$.

\proclaim{Theorem 6.1}The exponential map Exp: $\frak s \frak o (3,1)  
\rightarrow SO(3,1)^+$
given by $F \mapsto e^F$ is onto.  That is, for every proper Lorentz 
transformation $L$, there exists an $F  \in  \frak s \frak o (3,1)$ so that 
$L = e^F$.
\endproclaim

To prove the above theorem, we need to consider the complexification
 $\frak s \frak o (3,1) \otimes \Bbb C$ operating on
  $\Bbb R^{3,1} \otimes \Bbb C$. This last is isomorphic to $\Bbb C ^4$ and has an inner
product which is of the type $- + + +$ on $\Bbb R^{3,1}$ and extends to the complex 
vectors by $\langle i \vec v, \vec w\rangle = \langle \vec v, i \vec w\rangle = i \langle \vec v, \vec w\rangle$.
See \cite{Gottlieb(2001), Section 2} for more details.

Now let $c: \frak s \frak o (3,1) \rightarrow \frak s \frak o (3,1) \otimes \Bbb C$
given by $cF= F - i F^*$.  The image of $c$,
denoted $S$, is a three--dimensional complex vector space.  The set of 
operators of the form $aI + b cF$ will be denoted by $I + S$.  Note that $I + S$ is
a vector space isomorphic to $\Bbb R ^{3,1} \otimes \Bbb C$, and that $I + S$ is 
closed under multiplication, as the following lemma shows.

\proclaim {Lemma 6.2}Let $F$ and $G \in  S$ denote $cF$ and $cG$.  Then
$(a I + bF) (\alpha I + \beta G) = (a \alpha + b \beta \langle F,G\rangle) I + (b \alpha F + a \beta G + \ds{b \beta \over 2} [F,G])$
\endproclaim 
Now we say that $L \ \in \ I + S$ is a {\it biquaternion Lorentz transformation}
if $\langle Lu, Lv\rangle = \langle u, v\rangle$.  Any biquaternion Lorentz transformation $L$ must have the form $L= a I + bF$, where $F  \in  S$, 
such that $a^2 - b^2  \lambda^2_F =1$.

That is, $L^{-1} = a I - bF$.

\proclaim {Theorem 6.3}Every complex Lorentz transformation $L$ is an exponential, 
that is $L= e^F$ for some $F \in S$, except for $L = -I + N$ where 
$N  \in S $ is null, that is $N^2 = 0$.
\endproclaim 

\demo{Proof}Recall \cite{Gottlieb(1998), Theorem 8.5} where $F  \in  S$ 
that
$$
    e^F = \cosh (\lambda_F) I + \ds{\sinh (\lambda_F) \over \lambda_F} F
\tag {**}
$$
Now $L= aI + H$ where $H  \in  S$ and $a^2 - \lambda^2_H = 1$.  
So the first obstruction to showing that $L$ is an exponential is solving the 
equation
$\cosh (\lambda)= a$.  We shall show below that such a $\lambda$ always 
exists.  Next, if \ \ $\ds{\sinh (\lambda) \over \lambda} \not= 0$, then
$$
   L= aI + H = \cosh(\lambda) I + \ds{\sinh  \lambda \over \lambda} 
                \left(\ds{\lambda \over \sinh \lambda} H\right)
             =: \cosh (\lambda) I + \ds{\sinh \lambda \over \lambda} D = e^D
$$  
Hence $L$ may not be an exponential if \ \ $\ds{\sinh (\lambda) \over \lambda}=0$.

Now $\ds{\sinh \lambda \over \lambda} = 0$ exactly when $\lambda= \pi ni$ for
$n$ a non-zero integer. (Note that $\ds{\sinh (0) \over 0} =1$).  Then 
$$
   a = \cosh (\lambda) = \cosh (\pi ni)= \cos (\pi n) = (-1)^n.
$$
If $n$ is even, then $L= I + N = e^N$ where $N$ must be null.

If $n$ is odd, then $a= (-1)^n = -1$, so $L= -I + N$ where $N$ must be null
or zero.  Now $e^{B}= -I$ where $B  \in  S$ has eigenvalue 
$(2k + 1) \pi i$.  But $-I + N = -e^{-N}$ cannot be an exponential, because it has  a real 
eigenvector with negative eigenvalue.   This proves Theorem 6.3 except for the following lemma.
\enddemo

\proclaim{Lemma 6.4}
\item{a)}$\cosh(\lambda) = a$ always has a solution over the complex numbers.

\item{b)}$\sinh(\lambda) = 0$ if and only if $\lambda= \pi ni$.
\endproclaim

\demo{Proof}First we show b).  Now $\sinh (\lambda) = \ds{e^{\lambda} - 
e^{- \lambda} \over 2} = 0$.

Thus $e^{2 \lambda} = 1$, hence $2 \lambda = 2 \pi ni$ so $\lambda= \pi ni$.

Next we show a).
Now $\cosh(\lambda) = \ds{e^{\lambda} + e^{- \lambda} \over 2} = a$.  
Hence $(e^{\lambda})^2 - 2 ae^{\lambda} + 1 = 0$

Hence $e^{\lambda} = \ds{2a \pm \sqrt{4a^2 - 4} \over 2}= a \pm \sqrt{a^2 - 1}$.

Now $e^{\lambda} = b$ has a solution for all $b$ except $b=0$.  But $a \pm 
\sqrt{a^2 -1}$ cannot equal zero, hence we have shown there is a solution
for each $a$.
\enddemo

\demo{Proof of Theorem 6.1}We show the exponential map is onto $SO (3,1)^+$ 
by showing the products of two exponentials is an exponential.  That is 
$e^F e^G= e^D$ for $F, G, D \in \frak so (3,1)$.  Now $e^F = 
e^{{1 \over 2}cF} e^{{1 \over 2} \ov cF}$ where $\ov c F= F + i F^*$.  This follows
since $cF$ and $\ov c F$ commute.  Also for this reason, $e^{cF}$ and $ e^{\ov c G}$
commute.  Thus $e^F e^G = e^{{1 \over 2} cF} e^{{1 \over 2} cG} 
e^{{1 \over 2}\ov c F} e^{{1 \over 2} \ov c G}$.  Now $e^{{1 \over 2} cF} 
e^{{1 \over 2} cG}$ is a complex Lorentz transformation in $I+S$.  So either it 
is an exponential 
$e^{cD}$, or it has the form $-I + cN= - e^{cN}$ by Theorem 6.3.  Now Theorem
6.3 also holds for $I+\ov S $.  Hence we have 
$e^F e^G = e^{2D}$ or $e^F e^G= (-e^{cN}) (-e^{\ov c N})= e^{2N}$. \qed
\enddemo

\proclaim{Corollary 6.5} The exponential map $\text{Exp}:{\frak so}(3,1) \otimes \bC \to SO(\Bbb R^{3,1}\otimes\Bbb C)$ is not onto.  If $N\in \frak s \frak o(3,1)$ is null, then $-e^N$ is
not an exponential even though $-e^{cN}$ is an exponential.
\endproclaim

\demo{Proof} As explained in \cite{Gottlieb(2001)}, we can extend duality $F^*$ to
skew symmetric matrices $\pmatrix 0 & \vec E\\  \vec E & \times \vec B\endpmatrix$
where $\vec E$ and $\vec B$ are complex vectors.  Then $cF=F-iF^*$ and
$\ov cF=F+iF^*$ satisfy the same properties as in the complexification of the 
real case.  Now consider $e^Fe^G$ where $F$, $G\in S$.  Then $cF={1\over 2}\,cF+{1\over 2}\,\ov c F$, so $e^Fe^G=e^{{1\over 2}cF}
e^{{1\over 2}\ov cF}e^{{1\over 2}cG}e^{{1\over 2}\ov cG}$.  Now $cF=cA$ for some
$A\in \frak s \frak o(3,1)$, and $\ov c F= \ov cA'$ for $A' \in \frak s \frak o(3,1)$, hence
$$
e^Fe^G=e^{cA} e^{\ov cA'} e^{cB}e^{\ov cB'}=(e^{cA}e^{cB})(e^{\ov cA'}e^{\ov cB'}),
$$
and so $e^{cA}e^{cB}$ equals either $e^{cD}$ or $-e^{cN}$.  But $(-I)e^{cN}=e^{(2n+1)\pi i\ov cE}e^{cN}=e^{(2n+i)\pi i\ov c E+cN}$ where $E$ has eigenvalue
equal to $1$.  So in both cases $e^{cA}e^{cB}$ is an
exponential.

Now $-e^{cN}$ is an exponential since $-e^{cN}=e^{\pi i\ov cE}e^{cN}=
e^{\pi i\ov cE+cN}$ where $ E$ has eigenvalue $\lambda_{c E}=1$.
On the other hand $-e^N$, where $N$ is the real part of a null $cN$, cannot be an 
exponential, since if
$-e^N=e^F$, then  $s$, the unique eigenvector for $e^N$, applied to this
equation gives 
$-s=e^Fs=e^{\lambda_F}s$, so $\lambda_F=(2n+1)\pi i$ for some $n$.  Thus $F$
has another linear independent null eigenvector, which contradicts $-e^N$ having only one.
\qed \enddemo

\head 7. Eigenvectors
\endhead

In this section we give explicit formulas for the eigenvectors  and eigenvalues of proper
Lorentz transformations and their Lie algebra. We show the Doppler shift factor arises
as a kind of Berry's phase.

\proclaim{Theorem 7.1}  Let $F\in \frak s \frak o (3,1)$ and let $\lambda_F$ be an
eigenvalue of $F$ and $\lambda_T$ be an eigenvalue of $T_F$. The eigenvalue
of $cF$ is $\lambda_{cF} = \lambda _F -i \lambda_{F^*}$ and

\item{a)} $\ds\lambda_T=\sqrt{({E^2-B^2\over 2})^2+(\BE\cdot\BB)^2}$

\item{b)} $\ds\lambda_F=\pm\sqrt{\lambda_T+{(E^2-B^2)\over 2}}$,\quad
$\lambda_{F^*}=\pm\sqrt{\lambda_T-{(E^2-B^2)\over 2}}$.

\endproclaim
\demo{proof}
This is Theorem 5.4 of \cite{Gottlieb(1998)}.
\enddemo

Now the image of $\lambda_{cF} I + cF$ is the 2-dimensional space of eigenvectors of $cF$
with eigenvalue $\lambda_{cF}$. The image of 
$\lambda_{\overline cF} I + \overline cF$
is the 2-dimensional space of eigenvalues of $\overline cF$. Note that this is the complex 
conjugate of the eigenspace of $\lambda_{cF} I + cF$. Now let $u$ be a vector
of length $-1$ in the Minkowski metric, an observer in relativity theory.
Then $ s := (\lambda_{cF} I + cF)(\lambda_{\overline cF} I + \overline cF )u$
is in both eigenspaces, since the operators commute. And $s$ is a real vector since $u$ is. 
So $s$ is not only an 
eigenvector for $cF$ and $\overline cF$, but also for the real part
$F$ and the imaginary part $F^*$, and hence for the stress-energy tensor
$T_F$ and the Lorentz transformation $e^F$. See section 5, \cite{Gottlieb(1998)}.

\proclaim{Theorem 7.2} The eigenvector $s := (\lambda_{cF} I + cF)(\lambda_{\overline cF} I + \overline cF )u $ for $F\in \frak s \frak o (3,1)$ 
with 
$\BE = Fu$ and $\BB = -F^*u$
satisfies the
following equation:
$$
s=2\big((\lambda_T+{E^2+B^2\over 2}\ ) u+\BE\times\BB+\lambda_F\BE-\lambda_{F^*}\BB \big).
\tag7.1
$$
\endproclaim
\demo{proof}
This is Corollary 6.8 of \cite{Gottlieb(1998)}.
\enddemo

\proclaim{Corollary 7.3} 
For a null $N \in {\frak so}(3,1)$, the eigenvector is
$$s = 2 \big( ({E^2+B^2\over 2}\ ) u+\BE\times\BB \big) \tag7.2$$

\endproclaim

\demo{proof}
Now $N$ null is the real part of the null $cN$. So $ \lambda_{cN} = \lambda_N - i \lambda_{N^*} = 0$. Hence $\lambda_N = \lambda_{N^*} =  \lambda_T =0$. Then plug this into Theorem 7.2.
\qed
\enddemo

Since there are at most two eigenvalues $\lambda_{cF}$, one the negative of the other,
and since the null matrices have only one eigenvalue, 0, we see from  the above results that
there are two null real eigenvector spaces for the generic case and one null real eigenvector space for a null matrix.

Now we can use the above formulas to give us something like a connection on the 
eigenbundles of a field of $F \in \frak{so}(3,1)$ on Minkowski space-time. And we can consider
what occurs as we move around a closed time-like circuit in space-time, that is, two time-like
paths starting with the same velocity at time $0$ and ending at the same point at some positive
time. Then the eigenvectors formulas will progress according to the formulas until they meet
at a future time where they lie in the same 1-dimensional space, but they differ by a factor. We
can calculate that factor. It only depends upon the tangent velocities $u$ and $u'$ at the point
of intersection and the factor is real This differs from Berry's phase, in which the factor is complex
and usually depends upon the history of the paths, yet it has the same feel to it.

We follow Scholium 8.2 of \cite{Gottlieb(1998)}

Let $s_u$ be an eigenvector of $F$ corresponding to $\lambda_F$ as seen by an
observer $u$.
Suppose
$$
u'=\ds{1\over \sqrt{1-w^2}} (u+\bw)\tag7.3  
$$
is another observer.
Then $u'$ sees a different eigenvector $s_{u'}$.
But $s_{u'}$ must be a multiple of $s_u$ since they are eigenvectors.
So the question is, what is the multiple in terms of $\bE, \bB$ and $\bw$?
The answer is:

\proclaim{Theorem 7.4}
$$
s_{u'}={1\over \sqrt{1-w^2}} \left[ 1+
{ -(\bE\times \bB)\cdot \bw+\lambda_F \bE\cdot \bw-\lambda_{F^*}\bB\cdot \bw
 \over \lambda_T+\ds{E^2+B^2\over 2}} \right] s_u. \tag 7.4  
$$
\endproclaim

\demo{Proof} Define $$\varphi(v)={\langle v,s_-\rangle\over \langle u,s_-\rangle
} s_u\tag7.5  
$$
where $s_-$ is an eigenvector corresponding to $-\lambda_F$.
Then $\varphi$ is a linear map whose image is the span of $s_u$ and whose kernel
 is the space of vectors orthogonal to $s_-$.
Now $\varphi(u)=s_u$.

Now $\Phi:=(\lambda_{cF} I + cF)\circ (\overline{\lambda_{cF}} I + \overline{c}F
)$ has the same properties and let $\Phi (u):=s_u$.
Then $\Phi=\varphi$.
Let $s_-=\Phi_- (u)=(-\lambda_{cF} I+cF)\circ (\overline{-\lambda_{cF}} 
I + \overline{c}F)u$.

Now
$$
s_u=2\left(\lambda_Tu + {E^2+B^2\over 2}u+\bE\times \bB+\lambda_F \bE-\lambda_{F
^*} \bB\right)\tag7.6  
$$
from (7.2) and $s_-$ is the same with the signs changed on $\lambda_F$ and 
$\lambda_{F^*}$:

$$
s_-=2\left(\lambda_Tu + {E^2+B^2\over 2}u+\bE\times \bB-\lambda_F \bE+\lambda_{F
^*} \bB\right)\tag7.7  
$$

Now $s_{u'}=\varphi(u')=\ds{{\langle u',s_-\rangle\over 
\langle u,s_-\rangle} }s_u$.
Substituting (7.3) into this equation yields
$$
s_{u'}={1\over \sqrt{1-w^2}} \left(1+{\langle \bw,s_-\rangle\over 
\langle u,s_-\rangle}\right) s_u.\tag7.8 
$$
Now
$$
\langle u,s_-\rangle=-2 \left(\lambda_T+{E^2+B^2\over 2}\right)\tag7.9  
$$
using (7.7).
Then using (7.7) to calculate $\langle \bw,s_-\rangle$ and substituting this into (7.8) we obtain (7.4).\qed

\enddemo

Now (7.4) holds for all $F\in {\frak so}(3,1)$.
If we restrict to null $F$ we should see (7.4) reduce to a simpler form.
In the null case $\lambda_F=\lambda_{F^*}=0$ and $E=B$.
So equation (7.4) reduces to
$$
s_{u'}={1\over \sqrt{1-w^2}} \left(1-\bw\cdot { (\bE\times \bB)\over E^2}\right)
 s_u .\tag7.10  
$$
Now $\bw \cdot \ds{(\bE\times \bB)\over E^2}$ is the component along the 
$\bE\times \bB$ direction.
If we assume that $\bw=\bw_r$, that is $\bw$ is pointing in the radial direction, then
$$
s_{u'}=\sqrt{ 1-w_r\over 1+w_r} s_u.\tag7.11  
$$
Here $\ds{\sqrt{ 1-w_r\over 1+w_r}}$ is the Doppler shift ratio.
This suggests that null $F$ propagate along null geodesics by parallel 
translation.

Now the fact that $I + S$ and $I + \overline S$ commute leads to a richer situation in analogy
to Berry's phase considerations. If $V$ is a 2--dimensional eigenspace for $F \in I + S$, then it is invariant under any $G \in I+ \overline S$. In fact, any null 2--dimensional subspace of
complexified Minkowski space is either an eigenspace of an $F \in S$ or an eigenspace of an
$F \in \overline S$. The action of $x, y, z$ on $\overline V$ is an irreducible action of
the spin Lie algebra, and the action of $X, Y, Z$ on $V$ is also an irreducible action of the spin
Lie algebra on $V$. The particular basis of the actions have a sign difference which \cite{Ryder(1988)}
calls left and right spin 1/2 actions. 

Now, for example, the nullquat $(\lambda_{cF} I + cF)$ composed with
$e^{\overline cG}$ and applied to a vector $u$ must be an eigenvector of cF. So if these three
quantities are varied,  one gets a formula giving the progression of an eigenvector of $cF$.


\head 8 Physical examples of eigenvectors and quantum probability
\endhead

We will point out two examples of inner products of eigenvectors of
$F$ in $\frak M$ which give probabilities underlying two important cases
in \cite{Sudbery (1986)}: Page 200, equation (5.84) which gives the probability
of spin along an axis at angle $\theta$ from the spin direction of the 
particle. In this case the probability of spin $+ 1/2$ is equal to
the Minkowski innerproduct 
$$-{1\over 2} \langle  u + \bold v , u + \bold w \rangle = \sin^2(\theta/2)$$
where $u$ is an observer, i.e. $\langle u , u \rangle = -1$, and $\bold v $ and
$\bold w$ are unit vector in the rest space of $u$ pointing along the direction
of spin of the particle and the direction of the measurement, usually the
gradient of a pure $\bB$ field. Note both $ u + \bold v$ and $ u + \bold w$ are both 
null vectors, and hence possible eigenvectors of some operators in $\frak M$.

The other example is on P. 273, equation (6.121) of \cite{Sudbery (1986)}. Here the distribution
of electrons with specific velocity $v$ is given by $1 - v \cos(\theta)$,
where the electrons decay from a Cobalt $60$ atom in a strong magnetic field
$\bB$. Here $\theta$ is the angle between the magnetic field $B$ and the
velocity of the electron $v$. If we let $u$ represent the center of mass
observer $u$ and $u' = {1\over \sqrt{(1-v^2)}} (u + \bold v)$ represent the 4-velocity of the electron
and $u + {1 \over B} \bB$ be the normalised eigenvector of $F$ representing the 
pure $\bB$ field, then
$$ - \langle  \sqrt{(1-v^2)} u', u + {1 \over B} \bold B \rangle$$ 
equals this distribution.

\enddocument